# Piecewise smooth interval maps with nonvanishing derivative

Ale Jan Homburg
Institut für Mathematik I
Freie Universität Berlin
Arnimallee 2-6
14195 Berlin
Germany
alejan@math.fu-berlin.de

**Abstract**

We consider the dynamics of piecewise smooth interval maps $f$ with nowhere vanishing derivative. We show that if $f$ is not infinitely renormalizable, then all its periodic orbits of sufficiently high period are hyperbolic repelling. If in addition all periodic orbits of $f$ are hyperbolic, then $f$ has at most finitely many periodic attractors and there is a hyperbolic expansion outside the basins of these periodic attractors. In particular, if $f$ is not infinitely renormalizable and all its periodic orbits are hyperbolic repelling, then some iterate of $f$ is expanding. In this case, $f$ admits an absolutely continuous invariant probability measure.

## 1 Introduction

There is an extensive theory on the dynamics of one dimensional maps. Especially smooth maps have attracted much attention; the most detailed statements have been obtainerd for quadratic maps or for unimodal maps with negative Schwarzian derivative. See [MelStr,1993] for an account of the theory. The existence of critical points plays a dominant role in the dynamics of smooth maps. A natural question occurs what can be said on the dynamics of piecewise smooth maps for which the derivative nowhere vanishes. Studying the dynamics of such maps is the goal of this paper.

We study the dynamics of maps $f : I \to I$ on a compact interval $I$ with a finite number of turning points. A turning point is a local extremum in the interior of $I$. A continuous map $f$ is called multimodal if it possesses a finite number of turning points. Let $\mathfrak{E}$ be the class of multimodal maps $f : I \to I$, so that $f$ is strictly monotone outside the set of turning points $\mathfrak{T}$ and $\ln |Df|$ is Lipschitz continuous.

We show that from a metric point of view, the possible dynamics of a map $f \in \mathfrak{E}$ is limited. We show that if $f$ is not infinitely renormalizable then there exists $M > 0$ so that all periodic orbits of period larger than $M$ are hyperbolic repelling. To exclude pathological dynamics one can further assume that all periodic orbits are hyperbolic. Then there is at most a finite number of periodic attractors. We prove that in this situation there is an exponential expansion outside the basin of attraction of the periodic attractors.

Let us state our main result.

**Theorem 1.1** *Let $f \in \mathfrak{E}$ be at most finitely often renormalizable. Then there are numbers $K_n$ with $K_n \to \infty$ as $n \to \infty$, so that*

$$|Df^n(p_n)| \geq K_n \qquad (1)$$

*for each periodic point $p_n$ with minimal period $n$.*

*If moreover all periodic orbits are hyperbolic, then there is only a finite number of periodic attractors and there exist $C > 0, \lambda > 1$ so that for all $x \in I$ with $f^n(x)$ not in the immediate basin of attraction of an attracting periodic orbit,*

$$|Df^n(x)| \geq C\lambda^n.$$

*In particular, if $f$ only has periodic repellers, there exists $N > 0$ with*

$$|Df^N| > 1.$$

In section 5 we give variants of this theorem, for classes of maps with less smoothness requirements, but with additional restrictions on the dynamics.

Together with the results of Alsedà, López and Snoha [AlsLopSno,1995], classifying infinitely renormalizable piecewise smooth maps, this result provides a good understanding, from a metric point of view, of the dynamics of piecewise smooth multimodal maps with nowhere vanishing derivative.

As a corollary of the above theorem the following statement is obtained. A necessary and sufficient condition for $f \in \mathfrak{E}$ to be eventually expanding ($|Df^n| > 1$ for some $n \in \mathbb{N}$) is that $f$ is at most finitely often renormalizable and all its periodic orbits are hyperbolic repelling. It is well known that eventually expanding maps admit an absolutely continuous invariant probability measure (in short, an a.c.i.p.) [LasYor,1973], [MelStr,1993]. So the result in this paper shows that any piecewise smooth multimodal map with nowhere vanishing derivative admits an a.c.i.p. if all its periodic orbits are hyperbolic repelling and it is at most finitely often renormalizable. The property that some iterate of an interval map is expanding is persistent under smooth perturbations. Therefore, these a.c.i.p.'s occur persistently. This is in marked contrast with the situation for smooth maps, where a.c.i.p.'s do not occur persistently. We should remark though that in one parameter families of smooth multimodal maps satisfying some natural conditions, a.c.i.p.'s occur at a set of parameter values of positive measure [Jak,1981], [BenCar,1985], [MelStr,1993], [MarNow,1996], [Lyu,1996].

In the course of proving theorem 1.1 we present a simple proof of Mañé's theorem, giving an exponential expansion along orbits which stay outside a neighborhood of the set of turning points and do not converge to periodic attractors (see theorem 3.1). Actually, our proof of Mañé's theorem holds for a somewhat larger class of maps than the proofs of Mañé [Man,1985] and van Strien [Str,1990], [MelStr,1993] (compare also [Nus,1988]). Let $\mathfrak{C}$ be the class of maps on the interval $I$ defined as follows. A map $f$ is in $\mathfrak{C}$ if $f$ is $C^1$ except possibly at a finite set, there exists $C \geq 1$ so that $1/C \leq |Df| \leq C$, and $\ln |Df|$ has bounded variation. We prove Mañé's theorem for maps from $\mathfrak{C}$. This includes piecewise affine maps which have attracted some attention recently, see [GalMarTre,1994], [MarTre,1994], [LopSno,1995]. It is shown in [MarTre,1994], [LopSno,1995], that piecewise affine maps are not infinitely renormalizable. It was conjectured in [GalMarTre,1994] that piecewise affine maps with only hyperbolic repelling periodic orbits would be eventually expanding. Although some intermediate results in this paper, like Mañé's theorem, are proved for a class of maps including piecewise affine maps, our main result is proved for a class of maps that includes only simple piecewise affine maps. However, the statement of the main theorem also holds for maps $f \in \mathfrak{C}$ for which the limit sets of the turning points are not minimal Cantor sets, see section 5.

Let us say a few words on the proofs in this paper. A basic lemma we prove provides, for a map $f \in \mathfrak{C}$, a hyperbolic expansion outside the basin of the periodic attractors under the conditions that all periodic orbits are hyperbolic and there is a strong expansion along periodic orbits of high period. This lemma enables an easy proof of Mañé's theorem mentioned above. Indeed, it is readily seen that periodic orbits of high period staying outside a fixed neighborhood of the turning points have a strong expansion along them. Mañé's theorem follows from an application of the just described lemma.

From the above it is clear that a keyrole is played by the periodic orbits. Our strategy for proving theorem 1.1 is by demonstrating a strong expansion along periodic orbits of high period. It suffices to show this for periodic orbits that stay in the vicinity of the $\omega$-limit set of a turning point. Demonstrating this is fairly direct if the turning point is periodic or not recurrent. More work will be involved in establishing a strong expansion along periodic orbits of high period near the $\omega$-limit set of a nonperiodic but recurrent turning point. Dynamics near $\omega$-limit sets of recurrent turning points that are not minimal, is treated by direct arguments similar to the ones used in the previous sections. For periodic orbits near minimal $\omega$-limit sets we proceed as follows. If $\omega(c)$ is a Cantor set we prove, by adapting the existing proofs for $C^2$ maps, that it has zero Lebesgue measure. Using this we show that if $\omega(c)$ is a minimal Cantor set, it is a hyperbolic repelling set. It is then immediate that a strong expansion exists along periodic orbits of high period that stay near $\omega(c)$.

The organization of this paper is as follows. The next section contains some notation and tools that are used throughout the paper. In section 3 we prove the above mentioned theorem by Mañé. In section 4 we assume theorem 1.1 holds for orbits that stay near the $\omega$-limit sets of the turning points and show how to extend to all orbits. With the lemmas proved in section 4 one can easily treat multimodal maps for which turning points are not recurrent, or periodic, or have a hyperbolic repelling $\omega$-limit set. In section 5 we prove theorem 1.1 making use of the material in the previous sections, as well as results on the Lebesgue measure of $\omega$-limit sets. These results are collected in section 6. We show that the $\omega$-limit set of a turning point has zero Lebesgue measure if it contains no intervals. This result is basically due to [BloLyu,1990], [Var,1996], who treated $C^2$ multimodal maps.

Discussions I had with Henk Bruin, Gerhard Keller, Matthias St. Pierre, Duncan Sands (especially him) and Sebastian van Strien have been very enlightening. The book [MelStr,1993] by Welington de Melo and Sebastian van Strien was a valuable source of information. I thank the referee for his or her comments.

## 2 Prerequisites

This section collects some, mostly well known, properties of interval maps that we need in the sequel. We first define the three classes of maps occuring in this paper. Let $I$ be a compact interval.

Let $\mathfrak{E}$ be the class of multimodal maps $f : I \to I$, so that $f$ is strictly monotone outside the set of turning points $\mathfrak{T}$ and $\ln|Df|$ is Lipschitz continuous.

Let $\mathfrak{D}$ be the class of multimodal maps $f : I \to I$, so that $f$ is strictly monotone outside the set of turning points $\mathfrak{T}$ and $\ln|Df|$ restricted to each interval in $I\backslash\mathfrak{T}$ can be extended to a Lipschitz continuous map on a compact interval.

Let $\mathfrak{C}$ be the class of maps $f : I \to I$, so that $f$ is $C^1$ except possibly at a finite set, there exists $C \geq 1$ so that $1/C \leq |Df| \leq C$, and $\ln|Df|$ has bounded variation.

Note that $\mathfrak{E} \subset \mathfrak{D} \subset \mathfrak{C}$. In particular, for all maps $f$ in $\mathfrak{E}$, $\mathfrak{D}$ or $\mathfrak{C}$, $|Df|$ is bounded and bounded away from zero. For $f$ from $\mathfrak{C}$ or $\mathfrak{D}$, it is allowed that at a turning point $c$, $\lim_{x\downarrow c}|Df(x)|$ and $\lim_{x\uparrow c}|Df(x)|$ differ. This is not allowed for maps $f \in \mathfrak{E}$.

An interval $T \subset I$ is called a *homterval* if $f^i|_T$ is monotone for all $i \in \mathbb{N}$. $T$ is called a *wandering interval* if in addition $f^i(T) \cap f^j(T) = \emptyset$ for all $0 \leq i < j$. See [MelStr,1993] for the following lemma.

**Lemma 2.1 ('Contraction principle')** *Let $f : I \to I$ be a continuous multimodal map without wandering intervals. For any $\delta > 0$ there is a $\tilde{\delta} > 0$ so that for any interval $J \subset I$ with $|J| \geq \delta$ and for which $f^n(J)$ does not converge to some periodic orbit as $n \to \infty$, we have $|f^n(J)| \geq \tilde{\delta}$.* ∎

That we can apply the contraction principle follows from the following result from [MarMelStr,1992].

**Theorem 2.2** *A map $f \in \mathfrak{C}$ has no wandering intervals.* ∎

The derivative $Df$ might not exist in a point $x$. To avoid cumbersome notation, we will write e.g. $|Df^n(x)| \geq C$ for $\liminf_{y \to x} |Df^n(y)| \geq C$. The *distortion* of $f^n$ on an interval $J \subset I$ is defined as
$$\sup_{x,y \in J} |Df^n(x)| / |Df^n(y)|.$$
A collection $\mathfrak{I} = \{I_1, \ldots, I_L\}$ of subintervals of $I$ is said to have *intersection multiplicity* $S$ if the maximum over $x \in I$ of the number of intervals from $\mathfrak{I}$ containing $x$ is $S$.

**Lemma 2.3** *Let $f \in \mathfrak{C}$. For each $S > 0$ there is $D > 1$ so that for each interval $J \subset I$ with the collection $\{J, f(J), \ldots, f^{n-1}(J)\}$ having intersection multiplicity at most $S$, the distortion of $f^n$ on $J$ is bounded by $D$.*

PROOF. Denote by $K = \text{Var}(\ln |Df|)$ the variation of $\ln |Df|$. Let $\mathcal{L}$ be the set of points where $|Df|$ is discontinuous, let $L$ denote its cardinality. Write
$$M = \sup_{x \in \mathcal{L}} \left| \lim_{y \downarrow x} \ln |Df(y)| - \lim_{y \uparrow x} \ln |Df(y)| \right|.$$
For $x, y \in J$,
$$|\ln |Df^n(x)| - \ln |Df^n(y)|| \qquad (2)$$
$$= \left| \sum_{k=0}^{n-1} \ln |Df(f^k(x))| - \ln |Df(f^k(y))| \right|$$
$$\leq S(K + LM).$$

The distortion of $f^n$ on $J$ is thus bounded by $e^{S(K+LM)}$. ∎

Maps from $\mathfrak{D}$ have somewhat better distortion properties than maps from $\mathfrak{C}$. The Lipschitz constant of $f \in \mathfrak{D}$ is defined as the supremum of the Lipschitz constants of $f$ restricted to an interval of $I \setminus \mathfrak{T}$.

The proof of the following lemma goes just as the proof of lemma 2.3.

**Lemma 2.4** *Let $f \in \mathfrak{D}$. Denote the Lipschitz constant of $\ln |Df|$ by $K$. Suppose $J_n$ is an interval on which $|Df^n|$ is continuous. Then the distortion of $f^n$ on $J_n$ is bounded by $e^{K \sum_{i=0}^{n-1} |f^i(J_n)|}$.* ∎

Observe that the distortion of $f^n$ on $J_n$ with $J_n$ as in the above lemma, is close to 1 if $\sum_{i=0}^{n-1} |f^i(J_n)|$ is small. The following lemma tells how much we can extend an interval $J_n$ on which $f^n$ has bounded distortion, and still have bounded distortion on the larger interval.

**Lemma 2.5** *Let $f \in \mathfrak{D}$. There is a constant $K$ so that the following holds. Let $J_n = (a,b)$ and $T_n = (a,d) \supset J_n$ be intervals on which $|Df^n|$ is continuous. Let $\sigma = \sum_{i=0}^{n-1} |f^i(J_n)|$ and let $\tau > \sigma$ be such that $\sum_{i=0}^{n-1} |f^i(T_n)| = \tau$. Then*

$$\frac{|T_n \setminus J_n|}{|J_n|} \geq e^{-K\tau} \left( \frac{\tau}{\sigma} - 1 \right).$$

PROOF. Denote $\rho = |T_n \setminus J_n|/|J_n|$. By lemma 2.4, the distortion of $f^n$ on $T_n$ is bounded by $e^{K\tau}$ for some constant $K$. So, for $0 \leq m < n$ and some $z \in T_n \setminus J_n$,

$$\begin{aligned} |f^m(T_n)| &= |f^m(J_n)| + |Df^m(z)||T_n \setminus J_n| \\ &\leq |f^m(J_n)| + e^{K\tau} \frac{|f^m(J_n)|}{|J_n|} |T_n \setminus J_n| \\ &\leq |f^m(J_n)| + e^{K\tau} \frac{|f^m(J_n)|}{|J_n|} \rho |J_n|. \end{aligned}$$

Hence

$$\sum_{i=0}^{n-1} |f^i(T_n)| \leq \sum_{i=0}^{n-1} |f^i(J_n)| \left( 1 + e^{K\tau} \rho \right).$$

So $\rho$ satisfies $\tau/\sigma \leq 1 + e^{K\tau}\rho$, i.e. $\rho \geq e^{-K\tau}\left(\frac{\tau}{\sigma} - 1\right)$. ∎

A closed forward invariant subset $X \subset I$ of $f \in \mathfrak{C}$ is called *hyperbolic repelling* if there exists $C > 0, \lambda > 1$ with

$$|Df^n(x)| \geq C\lambda^n \qquad (3)$$

for all $x \in X$. A periodic point $p$ with minimal period $n$ is called hyperbolic if both $\lim_{y \uparrow p} |Df^n(y)|$ and $\lim_{y \downarrow p} |Df^n(y)|$ are either smaller or larger than one.

**Lemma 2.6** *Let $f \in \mathfrak{C}$. Let $X$ be a closed invariant set with the property that there exists $\lambda > 1$ so that for all $x \in X$, there exists $n_x$ with $|Df^{n_x}(x)| \geq \lambda$.*

*Then $X$ is hyperbolic repelling.*

PROOF. Let, for some $\tilde{\lambda}$, $1 < \tilde{\lambda} < \lambda$, $B_x \ni x$ denote the ball with $|Df^{n_x}|_{B_x}| > \tilde{\lambda}$. Because $X$ is a closed set, it can be covered by a finite number of balls $B_{x_1}, \ldots, B_{x_L}$. Write $N = \max\{n_{x_1}, \ldots, n_{x_L}\}$. This implies that for $x \in L$, $n \in \mathbb{N}$ there are positive integers $0 = n_0 < n_1 < \ldots < n_M \leq n$ with $n_j - n_{j-1}, n - n_M \leq N$ and

$$|Df^n(x)| = \left(\prod_{j=1}^{M} Df^{n_j - n_{j-1}}\left(f^{n_1 + \ldots + n_{j-1}}(x)\right)\right) Df^{n-n_M}\left(f^{n_1 + \ldots + n_M}(x)\right),$$

$$\tilde{\lambda} \leq |Df^{n_j - n_{j-1}}\left(f^{n_1 + \ldots + n_{j-1}}(x)\right)|.$$

With $N = \max\{n_{x_1}, \ldots, n_{x_L}\}$ and $m = \min_{x \in I} |Df(x)|$ we thus have

$$|Df^n(x)| \geq \tilde{\lambda}^M m^{n-n_M} \geq \tilde{\lambda}^{(n/N)} m^N.$$

Therefore there exists $N \in \mathbb{N}$ so that for each $x \in L$, $|Df^N(x)| \geq \lambda$. On a small neighborhood $\mathcal{V}$ of $L$, $|Df^N|_{\mathcal{V}}| \geq \tilde{\lambda}$. The lemma now easily follows. ∎

## 3 Mañé's theorem

For a subset $\mathcal{U}$ of $I$, write

$$\Gamma_n(\mathcal{U}) = \{x \in I, \; x, f(x), \ldots, f^n(x) \in I \backslash \mathcal{U}\}. \tag{4}$$

In the case of $C^2$ maps, the following theorem is due to Mañé [Man,1985], see also [MelStr,1993]. The argument given here seems a more direct one.

**Theorem 3.1** *Let $f \in \mathfrak{C}$. Let $\mathcal{U}$ be a neighborhood of the set of turning points $\mathfrak{T}$ of $f$. Then there are numbers $K_n$ with $K_n \to \infty$ as $n \to \infty$, so that*

$$|Df^n(p_n)| \geq K_n \tag{5}$$

*for each periodic point $p_n$ with minimal period $n$ and $\mathcal{O}(p_n) \subset I \backslash \mathcal{U}$.*

*If all periodic orbits in $I \backslash \mathcal{U}$ are hyperbolic, then there are $C > 0$, $\lambda > 1$, so that for each $n \in \mathbb{N}$ and $x \in \Gamma_n(\mathcal{U} \cup B_0)$,*

$$|Df^n(x)| \geq C \lambda^n. \tag{6}$$

*Here $B_0$ denotes the union of the immediate basins of the periodic attractors.*

PROOF. Let $p_n$ be a periodic point of $f$ with minimal period $n$ and $\mathcal{O}(p_n) \subset I\backslash\mathcal{U}$. If $n$ is large enough, $p_n$ will not be a periodic attractor with a turning point in its basin of attraction. Since $Df^n$ is the same at each point of $\mathcal{O}(p_n)$, we may replace $p_n$ by the point in $\mathcal{O}(p_n)$ closest to $\mathcal{U}$. We keep writing $p_n$ for this point. Let $J_n \ni p_n$ be the maximal interval with $f^n|_{J_n}$ monotone and $f^n(J_n) \cap \mathcal{O}(p_n) = \{p_n\}$. One easily sees that each point $x \in I$ is contained in at most two of the intervals $f^j(J_n)$, $0 \leq j < n$. By lemma 2.3, there exists $D > 0$ so that for all $n$, the distortion of $f^n$ on $J_n$ is bounded by $D$. Let $\delta$ be the minimal length of components of $\mathcal{U}$. If $\partial(f^i(J_n)) \cap \mathfrak{T} = \emptyset$ for all $i < n$, $f^n(J_n)$ contains a component of $\mathcal{U}$ and so $|f^n(J_n)| \geq \delta$. If $\partial(f^i(J_n)) \cap \mathfrak{T} \neq \emptyset$ for some $i < n$ with $i$ chosen the minimal number for which this holds, then $|f^i(J_n)| \geq \delta$. By the contraction principle, which we can apply since no turning point is in the basin of attraction of $\mathcal{O}(p_n)$, there exists a positive number $\tilde{\delta} < \delta$ so that $|f^n(J_n)| \geq \tilde{\delta}$. It follows that

$$|Df^n(p_n)| \geq \frac{1}{D}\frac{|f^n(J_n)|}{|J_n|} \geq \frac{1}{D}\frac{\tilde{\delta}}{|J_n|} \tag{7}$$

By theorem 2.2 $f$ has no wandering intervals. We can therefore apply lemma 3.2 below to get $|J_n| \to 0$ as $n \to \infty$; this proves (5).

By extending $f$ to a larger interval and altering $f$ in $\mathcal{U}$, we may assume that the turning points are in the basin of attraction of an attracting fixed point in $\partial I$. Then (5) holds, perhaps with different numbers $K_n$, for all periodic orbits. Assume now that all periodic orbits in $I\backslash\mathcal{U}$ are hyperbolic. From an application of lemma 3.3 below we obtain (6). ∎

**Lemma 3.2** *Let $f : I \to I$ be a continuous l-modal map. Suppose $f$ has no wandering intervals. For each $S > 0$, there exists $\xi_n > 0$ with $\xi_n \to 0$ as $n \to \infty$ so that for all intervals $J_n$ with $f^n|_{J_n}$ monotone, $f^n(J_n) \supset J_n$ and $\{J_n, \ldots, f_n(J_n)\}$ having intersection multiplicity bounded by $S$, we have $|J_n| \leq \xi_n$.*

PROOF. Assume, by contradiction, there exists $C > 0$ and a sequence $J_i$ with $|J_i| \geq C$, $f^i|_{J_i}$ is monotone, $f^i(J_i) \supset J_i$ and $\{J_i, \ldots, f^i(J_i)\}$ having intersection multiplicity bounded by $S$. Let $J$ be an interval contained in infinitely many $J_i$'s. So $J$ is a homterval. Since wandering intervals do not exist, $f^l(J) \cap f^k(J) \neq \emptyset$ for some $l < k \in \mathbb{N}$. We may choose $l, k$ minimal with this property. It is not hard to see that $l$ and $k$ are bounded by integers $l_0, k_0$ depending only on $|J|$. The interval $\bigcup_{s\in\mathbb{N}} f^{l+s(k-l)}(J)$ is a homterval, so any point in $J$ is attracted by a periodic point of minimal period $k - l$ or $2(k - l)$. Observe that $f^l(J_i)$ contains a periodic point $q_i$ of period at least $i/S$. If $v \in f^l(J)$, because $2(k - l) < i/S$ for $i$ large enough, there is $w \in (v, q_i)$, $0 \leq s < i$ with $f^s(W) \in T$. This contradicts that $f^i|_{J_i}$ is monotone. ∎

The following lemma gives a 'hyperbolic structure' outside the basins of periodic attractors if all periodic orbits are hyperbolic and there is a strong expansion along periodic orbits of high period.

**Lemma 3.3** *Let $f \in \mathfrak{C}$. Suppose there are numbers $K_n$ with $K_n \to \infty$ as $n \to \infty$, so that $|Df^n(p_n)| \geq K_n$ for each periodic point $p_n$ of minimal period $n$. If $\mathcal{U} \subset I$ is an open set and all periodic orbits in $I \backslash \mathcal{U}$ are hyperbolic, then there are $C > 0$, $\lambda > 1$, so that for each $n \in \mathbb{N}$ and $x \in \Gamma_n(B_0 \cup \mathcal{U})$, we have $|Df^n(x)| \geq C\lambda^n$. Here $B_0$ denotes the union of the immediate basins of the periodic attractors.*

PROOF. We show that

$$\limsup_{i \in \mathbb{N}} |Df^i(x)| = \infty \qquad (8)$$

for each point $x \in \Gamma_\infty(B_0 \cup \mathcal{U})$. Let us first finish the proof assuming (8). It follows from (8) that there exists $\lambda > 1$ so that for each $x \in \Gamma_\infty(B_0 \cup \mathcal{U})$, there is $n_x \in \mathbb{N}$ with $|Df^{n_x}(x)| > \lambda$. By lemma 2.6, there exists $C > 0$, $\lambda > 1$ so that $|Df^n(x)| \geq C\lambda^n$ for all $x \in \Gamma_\infty(B_0 \cup \mathcal{U})$. So there are $N > 0$ and a neighborhood $\mathcal{V}$ of $\Gamma_\infty(B_0 \cup \mathcal{U})$ so that $|Df^N|_\mathcal{V}| \geq \lambda$. We claim there exists $M \in \mathbb{N}$ so that $(I \backslash \mathcal{V}) \cap \Gamma_M(B_0 \cup \mathcal{U}) = \emptyset$. Indeed, there would otherwise be a sequence of points $x_i \to x$ in $I \backslash \mathcal{V}$ so that $f^j(x_i) \notin B_0 \cup \mathcal{U}$ for $0 \leq j \leq i$. But then $x \in \Gamma_\infty(B_0 \cup \mathcal{U})$, which is impossible. It follows that for some $N \in \mathbb{N}$ and $\tilde{\lambda} > 1$, $|Df^N| \geq \tilde{\lambda}$ on $\Gamma_N(B_0 \cup \mathcal{U})$. The lemma easily follows.

It remains to establish (8). Since (8) is clear for repelling periodic points, we may assume that $x$ is not periodic. If $f^l(x)$ is a turning point for some $l$, replace $x$ by $f^{l+1}(x)$. So we may assume $\mathcal{O}(x) \cap \{\mathfrak{T}\} = \emptyset$. Choose a point $y \in \omega(x)$ as follows. If $\omega(x) \cap \mathfrak{T} \neq \emptyset$, let $y$ be a turning point in $\omega(x)$. If $\omega(x) \cap \mathfrak{T} = \emptyset$, there exist $c \in \mathfrak{T}$ and an interval $(y, y')$ with $y \in \omega(x)$, $c \in (y, y')$, $f(y) = f(y')$ and $\omega(x) \cap (y, y') = \emptyset$. We can in fact choose $y$ and $y'$ so that $y' \notin \mathfrak{T}$. Indeed, since $\omega(x) \cap \mathfrak{T} = \emptyset$, altering $f$ near $\mathfrak{T}$ so that the values of $f$ at $\mathfrak{T}$ are slightly perturbed, doesn't change $\omega(x)$. Alter $f$ so that $f(\mathfrak{T}) \cap \omega(x) = \emptyset$ and no new transverse intersections of the graph of $f$ with $I \times \{y\}$ are created. Seek an interval $(y, y')$ as above for this altered map, this satisfies the required properties.

If $y$ is a turning point, define a function $\tau$ on a small neighborhood $\mathcal{W}$ of $y$ by $\tau(y) = y$ and $f(\tau(x)) = f(x)$ with $\tau(x) \neq x$ if $x \neq y$. If $y$ is not a turning point, let $\tau$ be a function defined on a small neighborhood $\mathcal{W}$ of $\{y\} \cup \{y'\}$ by $f(\tau(x)) = f(x)$ and $\tau(x) \neq x$. Note that $\tau(y) = y'$.

Two cases occur which have to be studied separately. Either there exist infinitely many points in $\mathcal{O}(x) \cap (y, \tau(y))$ or there are only finitely many points in $\mathcal{O}(x) \cap (y, \tau(y))$.

*There are only finitely many points in* $\mathcal{O}(x) \cap (y, \tau(y))$. Note that this case in particular occurs if $y \in \mathfrak{T}$, since then $(y, \tau(y)) = \emptyset$. Let $N$ be the maximal number so that $f^N(x) \in (y, \tau(y))$. Replace $x$ by $f^{N+1}(x)$.

The sets $V_k$ we now define will play an important role in establishing (8). For $k \in \mathbb{N}$ write

$$V_k = \{y \in \mathcal{W}, \quad f^i(y) \notin (y, \tau(y)), 0 < i < k, f^k(y) \in (y, \tau(y))\}. \tag{9}$$

We first discuss some properties of the sets $V_k$. Let $T_k$ be a connected component of $V_k$. We claim that the following two items hold.

- A point $a \in \partial T_k$ satisfies $f^k(a) \in \{a, \tau(a)\}$,
- $f^i(T_k) \cap f^j(T_k) = \emptyset$ for $0 \leq i < j < k$.

For the first item, observe that $a \in \partial T_k$ implies that for some $l, 0 < l \leq k$, $f^l(a) \in \{a, \tau(a)\}$. Take $l$ to be the minimal number with this property. Either $l = k$, then we are finished, or $l < k$. In the latter case, write $k = sl + t$ with $t < l$. Then $f^k(a) = f^{sl+t}(a) = f^t(a)$ or $f^t(\tau(a))$. So $t = 0$ by minimality of $l$ and $f^k(a) \in \{a, \tau(a)\}$. In particular, for each boundary point $a$ of $T_k$, either $a$ or $\tau(a)$ is periodic. For the second item, write $T_k = (a, b)$ with $b \in (a, \tau(a))$. If $f^i(T_k) \cap f^j(T_k) \neq \emptyset$ then there is $y \in T_k$ with $f^k(y) = f^{k-j+i}(a)$ or $f^k(y) = f^{k-j+i}(b)$. The first possibility contradicts the definition of $T_k$, the second possibility implies that either $T_k \cap f^{k-j+i}(T_k) \neq \emptyset$ or $\tau(T_k) \cap f^{k-j+i}(T_k) \neq \emptyset$. This gives $a \in f^{k-j+i}(T_k)$ resp. $\tau(a) \in f^{k-j+i}(T_k)$. Therefore there exists $z \in (a, b)$ with $f^{k-j+i}(z) = a$ resp. $f^{k-j+i}(z) = \tau(a)$. But then $f^k(z) = f^{j-i}(a)$, contradicting the definition of $T_k$.

Because the intervals $f^i(T_k)$, $0 \leq i < k$, are mutually disjoint, the minimal period of the periodic point $a$ or $\tau(a)$ with $a \in \partial T_k$, is $k/2$ or $k$.

Let $f^{n(i)}(x)$ be a sequence of closest returns to $(y, \tau(y))$; $n(0)$ satisfies $f^{n(0)}(x) \in \mathcal{W}$ and $n(i + 1)$ is the minimal integer so that $f^{n(i+1)}(x) \in (f^{n(i)}(x), \tau(f^{n(i)}(x)))$. Observe that $n(i+1) - n(i) \to \infty$ as $i \to \infty$ and $f^{n(i)}(x) \in V_{n(i+1)-n(i)}$. Let $T_{n(i+1)-n(i)}$ be the connected component of $V_{n(i+1)-n(i)}$ that contains $f^{n(i)}(x)$. Write

$$Df^{n(i+1)}(x) = \prod_{j=0}^{i} Df^{n(j+1)-n(j)}(f^{n(j)}(x)). \tag{10}$$

By lemma 2.3, for $z \in T_k$,

$$|Df^k(z)| \geq \frac{1}{D} C \min(K_{k/2}, K_k) \tag{11}$$

where $D$ bounds the distortion of $f^k$ on $T_k$ and $C = \min_{x \in \mathcal{W}} |Df(x)|/|Df(\tau(x))|$. Using this in (10) one gets $\limsup_{i \in \mathbb{N}} |Df^i(x)| = \infty$.

*There are infinitely many points in $\mathcal{O}(x) \cap (y, \tau(y))$*. Remains the case where $y \notin \mathfrak{T}$ and an infinite number of iterates of $x$ is contained in $(y, \tau(y))$. The reasoning is similar as above, but involves a different sequence of closest returns. Note that $\mathcal{O}(x) \cap (y, \tau(y))$ accumulates on $y$ or on $\tau(y)$. By replacing $y$ by $\tau(y)$ if necessary, we may assume that $\mathcal{O}(x) \cap (y, \tau(y))$ accumulates on $y$.

Define

$$V_k = \{z \in (y, c), \quad f^i(z) \notin (y, z), 0 < i < k, f^k(y) \in (y, z)\}. \tag{12}$$

We first study some properties of $V_k$. Let $T_k$ be a connected component of $V_k$. We claim that

- a point $a \in \partial T_k$ satisfies $f^k(a) = a$ or $f^k(a) = y$,
- $f^i(T_k) \cap f^j(T_k) = \emptyset$ for $0 \leq i < j < k$.

For the first item, let $a \in \partial T_k$. Then for some $l, 0 < l \leq k$, either $f^l(a) = a$ or $f^l(a) = y$. Let $l$ be the minimal number for which this holds. If $f^l(a) = a$ then, writing $k = sl + t$ with $t < l$, $f^k(a) = f^{sl+t}(a) = f^t(a)$ shows $t = 0$ and $f^k(a) = a$. Since $\mathcal{O}(y) \cap (y, \tau(y)) = \emptyset$ and $a \in (y, \tau(y))$, $f^l(a) = y$ and $f^k(a) \neq y$ for $l < k$ is not possible. So either $f^k(a) = a$ or $f^k(a) = y$. To obtain the second item, suppose by contradiction $f^i(T_k) \cap f^j(T_k) \neq \emptyset$. Then $f^k(T_k) \cap f^{k-j+i}(T_k) \neq \emptyset$. By minimality of $k$, $f^{k-j+i}(T_k)$ can not be contained in $f^k(T_k)$. So there exists $z \in T_k$ with $f^{k-j+i}(z) = a$, the boundary point of $T_k$ with largest distance to $y$. But then $f^k(z) = f^{j-i}(a)$ can not lie in $T_k$, a contradiction.

Let $f^{n(i)}(x)$ be a sequence of closest returns to $y$ in $(c, y)$; $n(0)$ satisfies $f^{n(0)}(x) \in (y, c) \cap \mathcal{W}$ and $n(i+1)$ is the minimal integer so that $f^{n(i+1)}(x) \in (y, f^{n(i)}(x))$. Note that $f^{n(i)}(x) \in V_{n(i+1)-n(i)}$. Let $T_{n(i+1)-n(i)}$ be the connected component of $V_{n(i+1)-n(i)}$ that contains $f^{n(i)}(x)$.

We claim that there exists a neighborhood $\mathcal{V}$ of $y$ so that $f^{n(i+1)-n(i)}$ is monotone on $T_{n(i+1)-n(i)}$ for all sufficiently large values of $i$. If such a neighborhood $\mathcal{V}$ would not exist, we could take a sequence of points $z_i \in \mathcal{V} \cap T_{n(i+1)-n(i)}$ converging to $y$ so that $f^{s_i}(z_i) = d$ for some $d \in \mathfrak{T}$, $0 \leq s_i \leq n(i+1) - n(i)$. Note $|z_i - f^{n(i)}(x)| \to 0$ and $|f^{n(i+1)-n(i)}(z_i) - f^{n(i+1)}(x)| \to 0$ as $i \to 0$. Because $\omega(x) \cap \mathfrak{T} = \emptyset$, there is $\delta > 0$ so that $|f^{s_i}(z_i) - f^{n(i)+s_i}(x)| = |d - f^{n(i)+s_i}(x)| \geq \delta$. This contradicts the contraction principle lemma 2.1, proving the claim. Using this claim it follows that for large $i$, $\partial T_{n(i+1)-n(i)}$ contains a periodic point. As above one concludes $\limsup_{i \in \mathbb{N}} |Df^i(x)| = \infty$. ∎

The use of the above lemma is not restricted to the proof of Mañé's theorem (i.e. to maps with the turning points contained in basins of periodic attractors). Also in the following sections our strategy will be to prove the existence of a strong expansion along periodic orbits of high period and then apply lemma 3.3. We will apply lemma 3.3 with $\mathcal{U}$ equal to the empty set. For clarity, let us formulate the corresponding lemma separately.

**Lemma 3.4** *Let $f \in \mathfrak{C}$. Suppose there are numbers $K_n$ with $K_n \to \infty$ as $n \to \infty$, so that $|Df^n(p_n)| \geq K_n$ for each periodic point $p_n$ of minimal period $n$. If moreover all periodic orbits are hyperbolic, then there are $C > 0$, $\lambda > 1$, so that for each $n \in \mathbb{N}$ and $x \in \Gamma_n(B_0)$, we have $|Df^n(x)| \geq C\lambda^n$. Here $B_0$ denotes the union of the immediate basins of the periodic attractors.* ∎

## 4 Hyperbolic limit sets

Mañé's theorem tells that our main theorem holds for orbits which stay outside a neighborhood of the turning points. In this section we show how to prove the main theorem if we assume it holds for orbits which stay near $\omega$-limit sets of turning points. This forms the contents of the following lemma. After statement and proof of this lemma we apply it to prove the main theorem in some simple situations, making assumptions on the orbits of the turning points.

**Lemma 4.1** *Let $f \in \mathfrak{C}$ and let $\mathfrak{T}$ denote the set of its turning points. Suppose that for each $c \in \mathfrak{T}$ there are a neighborhood $\mathcal{U}(\omega(c))$ of $\omega(c)$ and numbers $L_n$ with $L_n \to \infty$ as $n \to \infty$, so that*

$$|Df^n(p_n)| \geq L_n \tag{13}$$

*for each periodic point $p_n$ with minimal period $n$ and $\mathcal{O}(p_n) \subset \mathcal{U}(\omega(c))$. Then there are numbers $K_n$ with $K_n \to \infty$ as $n \to \infty$, so that*

$$|Df^n(p_n)| \geq K_n \tag{14}$$

*for each periodic point $p_n$ with minimal period $n$.*

PROOF. Write $\mathcal{U}$ for the union over the turning points of the neighborhoods $\mathcal{U}(\omega(c))$. Let $\mathcal{O}_n$ be a periodic orbit of minimal period $n$, containing a point in $I \backslash \mathcal{U}$. Since $Df^n$ is the same at each point of $\mathcal{O}_n$ we may, in order to prove (14) for $p_n \in \mathcal{O}_n$, replace $p_n$ by any point in $\mathcal{O}_n$.

We may assume that $p_n \in \mathcal{O}_n$ is contained in $I \backslash \mathcal{U}$. Write $J_n \ni p_n$ for the maximal interval with $f^n|_{J_n}$ monotone and $f^n(J_n) \cap \mathcal{O}_n = \{p_n\}$. Since the collection $\{f^i(J_n)\}$, $0 \leq iN,$, has intersection multiplicity at most 2, by lemma 2.3 $f^n$ has bounded distortion on $J_n$ with a bound $D$ not depending on $n$. Write $I_n \supset J_n$ for the maximal interval containing $p_n$ on which $f^n$ is monotone. Since $p_n \in I \backslash \mathcal{U}$ and $\partial f^n(I_n) \subset \mathcal{O}(\mathfrak{T})$, there is $\delta_0$ not depending on $n$ so that both components of $f^n(I_n)\setminus\{p_n\}$ have length at least $\delta_0$. By lemma 3.2, there is a sequence $\xi_n \to 0$, $n \to \infty$, with $|J_n| \leq \xi_n$. So (14) holds in case for some $p_n \in \mathcal{O}_n \cap I\setminus\mathcal{U}$, one point of $\partial J_n$ is contained in $\mathcal{O}(\mathfrak{T})$, since then $|f^n(J_n)| \geq \delta_0$.

It remains to prove (14) for periodic orbits $\mathcal{O}_n$ of minimal period $n$ so that for each $p_n \in \mathcal{O}_n \cap I\setminus\mathcal{U}$, $\partial f^n(J_n) \subset \mathcal{O}_n$. Suppose by contradiction that there is a constant $C > 0$ and a sequence of periodic points $p_{n(i)}$ in $I\setminus\mathcal{U}$ of minimal period $n(i)$, $n(i) \to \infty$ as $i \to \infty$, so that $|Df^{n(i)}(p_{n(i)})| \leq C$. From

$$|Df^{n(i)}(p_{n(i)})| \geq \frac{1}{D}|f^{n(i)}(J_{n(i)})|/|J_{n(i)}|$$

and lemma 3.2 it follows that $|f^{n(i)}(J_{n(i)})| \to 0$ as $i \to \infty$. Consider the set of limit points $\{f^j(p_{n(i)}), \ 0 \leq j < n(i), i \in \mathbb{N}\}$. Because $|f^{n(i)}(J_{n(i)})| \to 0$ as $i \to \infty$, this set of limit points contains an interval. We can therefore take a periodic point $y$ in $I\setminus\mathcal{U}$ contained in the interior of this set of limit points.

Denote by $k$ the minimal period of $y$. If $Df^k(y) < 0$, let $l = 2k$. Otherwise, let $l = k$. Let $P_1$ be a fundamental domain of $y$; $P_1$ is an interval of the form $[b, f^l(b))$ contained in the maximal interval around $y$ on which $f^l$ is monotone. Take $b$ to be an eventually periodic point. Let $P_n \subset (y, b)$ be such that $f^{n-1}(P_n) = P_1$. By lemma 4.2 below, there is a periodic point $q_{n(i)} \in \mathcal{O}(p_{n(i)})$, so that $|Df^j(q_{n(i)})| \leq C$ for all $j < n(i)$. Let $h(i)$ be the minimal integer so that $f^{h(i)}(q_{n(i)}) \subset P_n$. Let $H_n$ be the maximal interval containing $q_{n(i)}$ so that $f^{h(i)}$ is monotone on $H_n$ and $f^{h(i)}(H_n) \subset P_n$. We claim that for $n$ high enough,

- $f^{h(i)}(H_n) = P_n$,

- $f^q(H_n) \cap f^p(H_n) = \emptyset$,

for $0 \leq q < p < h(i)$. For the first item, if $f^{h(i)}(H_n) \neq P_n$, then $f^l(z) \in \mathfrak{T}$ for some $z \in H_n$ and $l < h(i)$. This contradicts $\mathcal{O}(\mathfrak{T}) \cap P_n = \emptyset$ for large $n$. For the second item, if $f^q(H_n) \cap f^p(H_n) \neq \emptyset$ for some $0 \leq q < p < h(i)$, then $f^{q+h(i)-p}(H_n) \cap f^{h(i)}(H_n) \neq \emptyset$. This contradicts $\mathcal{O}(\partial P_n) \cap P_n = \emptyset$ for $n$ high enough, which is a consequence of choosing $b \in \partial P_1$ eventually periodic.

From the above it follows that $f^{h(i)+n-1}(H_n) = P_1$ and $f^{h(i)+n-1}$ has uniformly bounded distortion on $H_n$. So, there exists $D > 0$ so that for all $i$,

$$|Df^{h(i)+n-1}(q_{n(i)})| \geq \frac{1}{D}\frac{|P_1|}{|H_n|}. \tag{15}$$

Since $|H_n| \to 0$ as $n \to \infty$, $|Df^{h(i)+n-1}(q_{n(i)})|$ is large if $i$ is large. This contradicts the definition of $q_{n(i)}$; since $h(i)+n-1$ is clearly bounded by $2n(i)$, we have that $|Df^{h(i)+n-1}(q_{n(i)})|$ is bounded by $C^2$. ∎

**Lemma 4.2** *Let $f \in \mathfrak{C}$. If $p_n$ is a periodic point of $f$ with $|Df^n(p_n)| \leq C$ for some $C > 1$, then there exists $q_n \in \mathcal{O}(p_n)$ with $|Df^j(q_n)| \leq C$ for all integers $j \leq n$.*

PROOF. Suppose by contradiction that there is a constant $\tilde{C} > C$ so that for all $x \in \mathcal{O}(p_n)$ there exists an integer $j(x) < n$ with $|Df^{j(x)}(x)| \geq \tilde{C}$. Let $x_1 = p_n$ and $x_2 = f^{j(x_1)}(x_1)$. Then $|Df^{n-j(x_1)}(x_2)| \leq C\backslash\tilde{C}$. Denote $x_3 = f^{j(x_2)}(x_2)$. Now either $j(x_1) + j(x_2) < n$ or $j(x_1) + j(x_2) > n$. In the first case, $|Df^{j(x_1)+j(x_2)}(x_1)| \geq \tilde{C}^2$. In the second case, $|Df^{j(x_1)+j(x_2)-n}(x_1)| \geq \tilde{C}^2\backslash C$. In both cases there is an integer $h(x_1) < n$ so that $|Df^{h(x_1)}(x_1)| \geq \hat{C}$ for some $\hat{C}$ which is at least a factor $\tilde{C}/C$ larger then $\tilde{C}$. Continuing this reasoning leads to a contradiction. ∎

In the following lemma we discuss expansion along periodic orbits near a hyperbolic repelling invariant set.

**Lemma 4.3** *Let $f \in \mathfrak{C}$ and let $c$ be a turning point of $f$. If $\omega(c)$ is a hyperbolic repelling set, then there are a neighborhood $\mathcal{U}$ of $\omega(c)$ and numbers $K_n$ with $K_n \to \infty$ as $i \to \infty$, so that for each periodic point $p_n$ of minimal period $n$ and with $\mathcal{O}(p_n) \subset \mathcal{U}$,*

$$|Df^n(p_n)| \geq K_n.$$

PROOF. We claim there exists $\delta_0 > 0$, $C > 0$ and $\lambda > 1$ so that for all $x \in I$ with $x, f(x), \ldots, f^n(x)$ contained in a $\delta_0$ neighborhood of $\omega(c)$, $|Df^i(x)| \geq C\lambda^i$, $0 \leq i \leq n$. Take $N$ so large that $|Df^N| > 3$ on $\omega(c)$. For each $x \in \omega(c)$ there is a ball $B_{\epsilon(x)}$ so that $|Df^N| > 2$ on $B_{\epsilon(x)}$. Since $\omega(c)$ is compact, it is covered by a finite set $\{B_{\epsilon(x_1)}, \ldots, B_{\epsilon(x_s)}\}$ of these balls. Let $\delta_0$ equal the minimum of $\epsilon(x_1), \ldots, \epsilon(x_s)$. Write $n = kN + l$ with $l < N$. Then $|Df^n(x)| \geq 2^k \min_{x \in I} |Df(x)|^l$. The claim and the proof of the lemma follow easily. ∎

Using the above lemma's one can easily treat dynamics of multimodal maps in $\mathfrak{C}$ for which the turning points are either periodic, or nonrecurrent, or have hyperbolic repelling limit set. This includes Misiurewicz maps; maps for which turning points are not recurrent. The dynamics of $C^2$ Misiurewicz maps was studied in [Str,1990], in that paper such maps were shown to admit an absolutely continuous invariant probability measure.

**Theorem 4.4** *Let $f \in \mathfrak{C}$. Suppose that for each turning point $c$ the following holds. Either $c$ is not recurrent, or $c$ is periodic, or $\omega(c)$ is a hyperbolic repelling set. Then there are numbers $K_n$ with $K_n \to \infty$ as $n \to \infty$, so that*

$$|Df^n(p_n)| \geq K_n \tag{16}$$

*for each periodic point $p_n$ with minimal period $n$.*

*If moreover all periodic orbits are hyperbolic, there is only a finite number of periodic attractors and there exist $C > 0, \lambda > 1$ so that for all $x \in I$ with $f^n(x)$ not in the immediate basin of attraction of an attracting periodic orbit,*

$$|Df^n(x)| \geq C\lambda^n. \tag{17}$$

*In particular, if $f$ only has periodic repellers, there exists $N > 0$ with*

$$|Df^N| > 1. \tag{18}$$

PROOF. We show that the assumption of lemma 4.1 is satisfied. The theorem then follows from lemma's 4.1 and 3.4.

Let $c$ be a turning point. If $\omega(c)$ is a hyperbolic repelling set, (13) follows from lemma 4.3. If $\omega(c)$ is periodic, periodic orbits which stay in a small neighborhood of it have the same period as $c$ or twice the period. So the estimate (13) holds for periodic points in a small neighborhood of $\omega(c)$ if $c$ is periodic. If $c$ is not recurrent, (13) holds for periodic orbits in a small neighborhood of $\omega(c)$ by Mañé's theorem 3.1. ∎

## 5 Proof of the main theorem

We make use in this section of a result on the measure of $\omega$-limit sets that we present in the next section. This result, an adaptation of work of Blokh and Lyubich [BloLyu,1989], [BloLyu,1990], Martens [Mar,1990], [Mar,1994] and Vargas [Var,1996], states that the $\omega$-limit set of any point of a map $f \in \mathfrak{D}$ either contains intervals or has zero Lebesgue measure. Using this result we derive a lemma stating that periodic orbits of high period near the $\omega$-limit set of a recurrent turning point, for a map $f \in \mathfrak{D}$ that is not infinitely renormalizable, are hyperbolic repelling. This is a first step in proving that such periodic points of high period are in fact strongly repelling.

Our main result is the following theorem. In its proof we make use of number of lemma's put after the proof.

**Theorem 5.1** *Let $f : I \to I$ be a map satisfying one of the following presumptions.*

1. *$f \in \mathfrak{E}$ is at most finitely often renormalizable,*

2. *$f \in \mathfrak{D}$ is unimodal and at most finitely often renormalizable,*

3. *$f \in \mathfrak{C}$ is such that the $\omega$-limit sets of the turning points are not minimal Cantor sets.*

*Then there are numbers $K_n$ with $K_n \to \infty$ as $n \to \infty$, so that*

$$|Df^n(p_n)| \geq K_n \qquad (19)$$

*for each periodic point $p_n$ with minimal period $n$.*

*If moreover all periodic orbits are hyperbolic, there is only a finite number of periodic attractors and there exist $C > 0, \lambda > 1$ so that for all $x \in I$ with $f^n(x)$ not in the immediate basin of attraction of an attracting periodic orbit,*

$$|Df^n(x)| \geq C\lambda^n.$$

*In particular, if $f$ only has periodic repellers, there exists $N > 0$ with*

$$|Df^N| > 1.$$

PROOF. It suffices to prove the theorem for maps that are not renormalizable. We prove that for each turning point $c$ of $f$ there is a neighborhood $\mathcal{U}$ of $c$ so that

$$|Df^n(p_n)| \geq K_n \qquad (20)$$

for periodic points $p_n$ of minimal period $n$ whose orbits are contained in $\mathcal{U}$. The theorem then follows from lemma 4.1.

Let $c \in \mathfrak{T}$. We have already seen that (20) holds if $c$ is not recurrent or periodic, compare theorem 4.4. Assume now that $c$ is recurrent and not periodic. The case where $\omega(c)$ is a minimal Cantor set is treated separately.

$\omega(c)$ *is a minimal Cantor set.* We prove that

$$\limsup_{n \in \mathbb{N}} |Df^n(x)| = \infty \qquad (21)$$

for all $x \in \omega(c)$. Since $\omega(c)$ is minimal, for each $x \in \omega(c)$ we have $c \in \omega(x)$.

Let $\tau$ be a function on a neighborhood $\mathcal{V}$ of $c$, defined by $\tau(c) = c$ and $f(\tau(y)) = f(y)$ with $\tau(y) \neq y$ if $y \neq c$. For $x \in \omega(c)$, let $f^{n(i)}(x)$ be a sequence of closest returns to $c$; $n(0)$ is such that $f^{n(0)}(x) \in \mathcal{V}$ and $n(i+1)$ is the smallest integer larger then $n(i)$ with $f^{n(i+1)}(x) \in (f^{n(i)}(x), \tau(f^{n(i)}(x)))$.

Write

$$|Df^{n(i+1)}(x)| = \prod_{j=0}^{i} Df^{n(j+1)-n(j)}(f^{n(j)}(x)). \qquad (22)$$

Let $V_k$ be defined as in (9);

$$V_k = \{y \in \mathcal{V}, \quad f^i(y) \notin (y, \tau(y)), 0 < i < k, f^k(y) \in (y, \tau(y))\}.$$

Here $\tau$ is a function on a neighborhood $\mathcal{V}$ of $c$, defined by $\tau(c) = c$ and $f(\tau(y)) = f(y)$ with $\tau(y) \neq y$ if $y \neq c$. Let $T_{n(i+1)-n(i)}$ be the component of $V_{n(i+1)-n(i)}$ containing $f^{n(i)}(x)$. It was shown in the proof of lemma 3.4 that

- $f^i(T_{n(i+1)-n(i)}) \cap f^j(T_{n(i+1)-n(i)}) = \emptyset$ for $0 \leq i < j < n(i+1) - n(i)$.

By the contraction principle (lemma 2.1),

$$\sup_{0 \leq j \leq n(i+1)-n(i)} |f^j(T_{n(i+1)-n(i)})| \to 0 \tag{23}$$

as $i \to \infty$.

By theorem 6.1, the Lebesgue measure of $\omega(c)$ is 0. Hence, the Lebesgue measure $|\mathcal{U}_\epsilon|$ of the $\epsilon$-neighborhood $\mathcal{U}_\epsilon$ of $\omega(c)$ satisfies

$$\lim_{\epsilon \to 0} |\mathcal{U}_\epsilon| = 0. \tag{24}$$

Because the intervals $f^l(T_{n(i+1)-n(i)})$ are mutually disjoint for $0 \leq l < n(i+1) - n(i)$ we conclude from (24) and (23) that $\sum_{l=0}^{n(i+1)-n(i)} |f^l(T_{n(i+1)-n(i)})| \to 0$ as $i \to \infty$. This implies that the distortion of $f^{n(i+1)-n(i)}$ on monotone branches of $T_{n(i+1)-n(i)}$ goes to 1 as $i \to \infty$.

If $f$ is a multimodal map from $\mathfrak{D}$, either $\partial T_{n(i+1)-n(i)}$ or $\tau(\partial T_{n(i+1)-n(i)})$ contains a periodic point. From lemma 5.2, using continuity of $|Df|$, it follows that the terms in the product (22) are strictly larger than 1 for $j$ large; (21) follows. If $f$ is a unimodal map from $\mathfrak{E}$, $f^{n(i+1)-n(i)}$ is monotone on $T_{n(i+1)-n(i)}$ and therefore $\partial T_{n(i+1)-n(i)}$ contains a periodic point. Hence (21) follows from lemma 5.2.

By lemma 2.6, $\omega(c)$ is a hyperbolic set. An application of lemma 4.3 yields (20).

*c is recurrent, nonperiodic, and $\omega(c)$ is not a minimal set.* The following reasoning applies to $f \in \mathfrak{C}$.

Let $c_1, \ldots, c_s$ with $c_s = c$ be the turning points in $\omega(c)$ with $\omega(c_i) \subset \omega(c_{i+1})$, $1 \leq i < s$. Note that there can still be other turning points $d \in \omega(c)$ with $\omega(d)$ strictly contained in $\omega(c)$.

Assume, by contradiction, the existence of $C > 0$ and a sequence of periodic orbits $\mathcal{O}_{n(i)}$ with minimal period $n(i)$ and with the maximal distance between $\mathcal{O}_{n(i)}$ and $\omega(c)$ converging to 0 as $i \to \infty$, so that $|Df^{n(i)}(p_{n(i)})| \leq C$ for $p_{n(i)} \in \mathcal{O}_{n(i)}$. By lemma 4.2, for each $i$ there exists a point $q_{n(i)} \subset \mathcal{O}_{n(i)}$ so that $|Df^j(q_{n(i)})| \leq C$ for all integers $j \leq n(i)$.

Let neighborhoods $S_n(c_h)$ of $c_h$ be as in lemma 5.3. Employing an induction argument we may assume that (20) holds for periodic orbits near $\omega$-limit sets of turning points that are strictly contained in $\omega(c_1)$. From this and Mañé's theorem 3.1, the minimal distance between $\mathcal{O}_{n(i)}$ and $\{c_1, \ldots, c_s\}$ goes to 0 as $i \to \infty$. So for any $n$ and for any $c_h$ there exists $i$ with $\mathcal{O}_{n(i)} \cap S_n(c_h) \ne \emptyset$. Let $h(i)$ be the minimal integer so that for some $j$, $f^{h(i)}(q_{n(i)}) \in S_n(c_j)$. Let $H_n$ be the maximal interval containing $q_{n(i)}$ on which $f^{h(i)}$ is monotone and with $f^{h(i)}(H_n) \subset S_n(c_j)$. We claim that

- $f^{h(i)}(H_n) = S_n(c_j)$,

- $f^k(H_n) \cap f^l(H_n) = \emptyset$ for $0 \le k < l < h(i)$.

Suppose the first item were false. Then $f^j(H_n)$ would intersect $\mathfrak{T}$ for some $j < h(i)$. By the contraction principle, $f^j(H_n)$ would then in fact hit a turning point $c_k$; the orbits of other turning points do not come near any $S_n(c_j)$. Since $f^j(q_{n(i)}) \notin S_n(c_k)$, there is $z \in H_n$ with $f^j(z) \in \partial S_n(c_k)$. This contradicts that $f^l(\partial S_n(c_k)) \cap S_n(c_j) = \emptyset$ unless $f^l(\partial S_n(c_k)) \subset S_n(c_j)$, for all positive integers $l$. To establish the second item, suppose $f^k(H_n) \cap f^l(H_n) \ne \emptyset$ for some $0 \le k < l < h(i)$. Then $f^{h(i)-l+k}(H_n) \cap S_n(c_k) \ne \emptyset$. By minimality of $h(i)$, $f^{h(i)-l+k}(H_n)$ can not be contained in $S_n(c_k)$, so that some $z \in H_n$ is mapped into $\partial S_n(c_k)$ by $f^{h(i)-l+k}$. A contradiction is derived as above.

By lemma 2.3, $f^{h(i)}$ has uniformly bounded distortion on $H_n$. From this and lemma 5.3, there is a constant $D > 0$ so that for all $i$,

$$|Df^{h(i)+s_j(n)}(q_{n(i)})| \ge \frac{1}{D} \frac{|f^{s_j(n)}(S_n(c_j))|}{|H_n|} \ge \frac{1}{D} \frac{\delta}{|H_n|}. \tag{25}$$

Since there are no homtervals, $|Q_n| \to 0$ as $n \to \infty$. Hence, $|Df^{h(i)+s_j(n)}(q_{n(i)})|$ is large if $n$ is large. From lemma 5.3 we obtain that the orbit piece $\{q_{n(i)}, \ldots, f^{h(i)+s_j(n)}(q_{n(i)})\}$ hits the interval $S_n(c_j)$ only once. Therefore $h(i) + s_j(n) \le 2n(i)$, from which it follows that $|Df^{h(i)+s_j(n)}(q_{n(i)})|$ is bounded by $C^2$, contradiction. ∎

**Lemma 5.2** *Let $f \in \mathfrak{D}$ be at most finitely often renormalizable. Suppose $c$ is a turning point of $f$ so that $\omega(c)$ is a Cantor set. Then there exist $\lambda > 1$, $\varepsilon > 0$ and $N \in \mathbb{N}$ so that for all periodic points $p_n$ of minimal period $n > N$ with $\mathcal{O}(p_n) \subset \mathcal{U}_\varepsilon$, where $\mathcal{U}_\varepsilon$ is the $\varepsilon$ neighborhood of $\omega(c)$,*

$$|Df^n(p_n)| \ge \lambda. \tag{26}$$

PROOF. We may assume that $f$ is not renormalizable.

Suppose a sequence of periodic points $p_{n(i)}$ of minimal period $n(i)$ accumulating on $\omega(c)$ exists so that

$$|Df^{n(i)}(p_{n(i)})| = \lambda_i \tag{27}$$

with $\lambda_i \to 1$ as $i \to \infty$. By Mañé's theorem 3.1 the minimal distance between $\mathcal{O}(p_{n(i)})$ and the set of turning points $\mathfrak{T}$ goes to 0 as $i \to \infty$. By taking a subsequence, we may assume that $p_{n(i)}$ converges to a turning point $\tilde{c} \in \omega(c)$. Let $\tau$ be a function defined on a small neighborhood $\mathcal{U}$ of $\tilde{c}$ so that $\tau(\tilde{c}) = \tilde{c}$ and $f(\tau(y)) = f(y)$ with $\tau(y) \neq y$ if $y \neq \tilde{c}$. We may assume that $p_{n(i)}$ is such that $(p_{n(i)}, \tau(p_{n(i)})) \cap \mathcal{O}(p_{n(i)}) = \emptyset$. Let $m(i) = n(i)$ if $Df^{n(i)}(p_{n(i)}) > 0$ and $m(i) = 2n(i)$ if $Df^{n(i)}(p_{n(i)}) < 0$. So $Df^{m(i)}(p_{n(i)}) > 0$.

Write $S_{m(i)}$ for the maximal interval in $(p_{n(i)}, \tau(p_{n(i)}))$ bounded by $p_{n(i)}$ with $f^{m(i)}(S_{m(i)}) \subset (p_{n(i)}, \tau(p_{n(i)}))$. We claim that

- $f^l(S_{m(i)}) \cap f^k(S_{m(i)}) = \emptyset$ for $0 \leq l < k < m(i)$.

Indeed, if $f^l(S_{m(i)}) \cap f^k(S_{m(i)}) \neq \emptyset$ for some $0 \leq l < k < m(i)$, then $f^{m(i)-k+l}(S_{m(i)}) \cap f^{m(i)}(S_{m(i)}) \neq \emptyset$. Since $\mathcal{O}(p_{n(i)}) \cap (p_{n(i)}, \tau(p_{n(i)})) = \emptyset$ there is $a \in S_{m(i)}$ with $f^{m(i)-k+l}(a) \in \{p_{n(i)}, \tau(p_{n(i)})\}$. But then $f^{m(i)}(a)$ can not lie in $S_{m(i)}$.

Write $H_{m(i)}$ for the maximal interval in $S_{m(i)}$ bounded by $p_{n(i)}$ on which $f^{m(i)}$ is monotone. Observe that

$$f^{m(i)}(H_{m(i)}) = f^{m(i)}(S_{m(i)}). \tag{28}$$

By theorem 6.1, the Lebesgue measure of $\omega(c)$ is 0. Hence, the Lebesgue measure $|\mathcal{U}_\epsilon|$ of the $\epsilon$-neighborhood $\mathcal{U}_\epsilon$ of $\omega(c)$ satisfies

$$\lim_{\epsilon \to 0} |\mathcal{U}_\epsilon| = 0. \tag{29}$$

By the contraction principle (lemma 2.1),

$$\sup_{0 \leq j \leq m(i)} |f^j(S_{m(i)})| \to 0 \tag{30}$$

as $i \to \infty$. Because the intervals $f^l(S_{m(i)})$ are mutually disjoint for $0 \leq l < m(i)$ we conclude from (29) and (30) that $\sum_{l=0}^{m(i)-1} |f^l(S_{m(i)})| \to 0$ as $i \to \infty$. This implies that the distortion of $f^{m(i)}$ on $H_{m(i)}$ goes to 1 as $i \to \infty$, see lemma 2.4. Because $Df^{m(i)}(p_{n(i)}) \to 1$ as $i \to \infty$, it follows that $|Df^{m(i)}| \to 1$ uniformly on $H_{m(i)}$ as $i \to \infty$. It is easily seen from this and (28) that $f$ is renormalizable. ∎

**Lemma 5.3** *Let $f \in \mathfrak{C}$ and let $c_1$ be a nonperiodic but recurrent turning point of $f$ so that $\omega(c_1)$ is not a minimal set. Let $c_1, \ldots, c_s$ be the turning points with $\omega(c_1) = \ldots = \omega(c_s)$. Then there exist $\delta > 0$ so that for each turning point $c_i$, there are decreasing neighborhoods $S_n(c_i)$ of $c_i$ and integers $s_i(n)$ such that*

$$|f^{s_i(n)}(S_n(c_i))| \geq \delta \tag{31}$$

*and $f^{s_i(n)}$ has uniformly bounded distortion on $S_n(c_i)$.*

*Furthermore, for n sufficiently large,*

$$f^l(S_n(c_i)) \cap f^k(S_n(c_i)) = \emptyset \tag{32}$$

*for all $0 \leq l < k < s_i(n)$ and*

$$f^l(\partial S_n(c_i)) \cap S_n(c_j) = \emptyset \text{ or } f^l(\partial S_n(c_i)) \subset S_n(c_j) \tag{33}$$

*for all $c_i \neq c_j$, $l \in \mathbb{N}$.*

PROOF. If $\omega(c_1)$ contains an interval, periodic orbits lie dense in it. If $\omega(c_1)$ is not a minimal Cantor set, it contains a a minimal set. Such a minimal set is a periodic orbit or a minimal Cantor set.

**(i)** $\omega(c_1)$ is a minimal Cantor set that contains a periodic point.

Let $y$ be a periodic point in $\omega(c_1)$. Write $k$ for the minimal period of $y$. If $Df^k(y) < 0$, let $l = 2k$, otherwise let $l = k$. Let $z_1$ be an eventually periodic point so that $f^l$ is monotone on $(y, z_1)$. Let $p_1 \in (y, z_1)$ be such that $f^l(p_1) = z_1$. Write $P_1 = [p_1, z_1)$ and let $P_n \subset (y, z_1)$ be such that $f^{(n-1)l}(P_n) = P_1$. Note that $P_n$ is a fundamental domain for $y$; $f^l(\text{int } P_n) \cap \text{int } P_n = \emptyset$ and $f^l(\partial P_n) \cap \partial P_n \neq \emptyset$.

We first construct $S_n(c_i)$ for a single turning point $c_i$. For positive integers $n$, let $\sigma(n)$ be the minimal integer with $f^{\sigma(n)}(c_i) \in P_n$. Let $\Sigma_n$ be the maximal interval containing $c_i$ so that $f^{\sigma(n)}(\Sigma_n) \subset P_n$. We claim that for $n$ sufficiently large,

- $f^i(\Sigma_n) \cap f^j(\Sigma_n) = \emptyset$ for $0 \leq i < j < \sigma(n)$,

Suppose by contradiction that $f^i(\Sigma_n) \cap f^j(\Sigma_n) \neq \emptyset$. Then $f^{\sigma(n)}(\Sigma_n) \cap f^{\sigma(n)-j+i}(\Sigma_n) \neq \emptyset$. Because $f^{\sigma(n)}(\Sigma_n) \subset P_n$ and $f^{\sigma(n)-j+i}(c_i) \notin P_n$, some $z \in \Sigma_n$ is mapped by $f^{\sigma(n)-j+i}$ to a point in $\partial P_{n+1} \cup P_n$. This contradicts $\mathcal{O}(z_1) \cap P_n = \emptyset$ for $n$ high which follows from the fact that $z_1$ is eventually periodic.

We may take $z_1$ so that $\omega(z_1)$ is outside $\omega(c_1)$. Then there is a neighborhood $T_1$ of $\partial P_1$ so that $\mathcal{O}(c_1) \cap T_1 = \emptyset$. Therefore, $f^{\sigma(n)}(\Sigma_n)$ contains a connected component of $T_1 \cap P_1$. By choosing $s_j(n) = \sigma(n)$ and $S_n(c_j) = \Sigma_n$, (31) is satisfied.

Note that $f^l(S_n(c_i))$ can not contain a boundary point of $S_n(c_j)$ for $0 \leq l < s_i(n)$ and $i \neq j$, since these boundary points are never mapped into $P_n$. By maximality of the intervals $S_n(c_i)$ it follows that, for $0 \leq l < s_i(n)$ and $i \neq j$, either $f^l(S_n(c_i)) \cap S_n(c_j) = \emptyset$, or $f^l(S_n(c_i)) = S_n(c_j)$; (33) is an easy consequence.

**(ii)** $\omega(c_1)$ is a Cantor set that contains no periodic points.

The construction is very much the same as the case where $\omega(c_1)$ is a Cantor set containing a periodic orbit, only different intervals $P_n$ are chosen. Let $y \in \omega(c_1)$ be a point whose $\omega$-limit set is minimal. As above we define a sequence of intervals $P_n$ with $|P_n| \to 0$ as $n \to \infty$, so that $f^n(P_n)$ is a fixed interval $P_1$. Choose an interval $P_1$ with the following properties. Let $\mathcal{U}$ be a small neighborhood of $\omega(y)$. Choose $P_1$ so that $P_1 \subset \mathcal{U}$, $P_1 \cap \omega(y) = \emptyset$, $P_1 \cap \omega(c_1) \neq \emptyset$ and $\partial P_1$ consists of eventually periodic points so that $\mathcal{O}(\partial P_1) \cap P_1 = \emptyset$. It is not hard to see that such an interval $P_1$ can be chosen. For $n > 1$, let $P_n$ be a maximal interval in $f^{-1}(P_{n-1})$ that is contained in $\mathcal{U}$ and intersects $\omega(y)$, so that $f^n$ is monotone on $P_n$. Because $\omega(y)$ is a hyperbolic repelling set, the intersection multiplicity of the collection of intervals $f^i(P_n)$, $0 \leq i < n-1$, is bounded.

It follows from lemma 2.3 that $f^n$ has uniformly bounded distortion on $P_n$. If $T_1$ is a neigborhood of $\partial P_1$ with $T_1 \cap \omega(y) = \emptyset$, then $f^n(P_n)$ contains a connected component of $P_1 \cap T_1$. We can conclude from the fact that $\omega(y)$ is a hyperbolic repelling set, that $|P_n| \to 0$ as $n \to \infty$. The rest is as above.

(iii) $\omega(c_1)$ contains an interval.

If $\omega(c_1)$ contains an interval we have a bit more work to do. Let $y$ be a periodic point in $\omega(c_1)$. Intervals $P_1 = [p_1, z_1)$ with $z_1$ eventually periodic and $P_n$ with $f^{(n-1)l}(P_n) = P_1$ are defined as in the case that $\omega(c_1)$ is a Cantor set containing a periodic point. For positive integers $n$, let $\sigma(n)$ be the minimal integer with $f^{\sigma(n)}(c_i) \in P_n$. Let $\Sigma_n$ be the maximal interval containing $c_i$ so that $f^{\sigma(n)}(\Sigma_n) \subset P_{n+1} \cup P_n$. As before one shows that for $n$ sufficiently large,

- $f^i(\Sigma_n) \cap f^j(\Sigma_n) = \emptyset$ for $0 \leq i < j < \sigma(n)$.

Let the integer $m_1$ be so that $f^{m_1}(z_1) = q_1$ is periodic. Write $n_1$ for the minimal period of $q_1$. If $Df^{n_1}(q_1) < 0$, we replace $n_1$ by $2n_1$, so that always $Df^{n_1}(q_1) > 0$. Let $b \in P_1$ be close enough to $z_1$ so that $f^{n_1}\big|_{(f^{m_1}(b), q_1)}$ is monotone.

*Either $P_2$ or $(b, z_1)$ is contained in $f^{\sigma(n)+(n-1)l}(\Sigma_n)$.* Let $S_n(c_i) = \Sigma_n$ and $s_i(n) = \sigma(n) + (n-1)l$.

*$f^{\sigma(n)+(n-1)l}(\Sigma_n)$ is contained in $(b, z_1)$.* Note that $f^{\sigma(n)+(n-1)l}(c_i) \in (b, z_1)$. Take a fundamental domain $P_1' \subset f^{m_1}(b, z_1)$ for $q_1$. Let $P_m'$ be the interval in $f^{m_1}(b, z_1)$ so that $f^{(m-1)n_1}(P_m') = P_1'$. Note that $f^{\sigma(n)+(n-1)l+m_1}(c_i)$ is contained in some interval $P_m'$. Let $\Sigma_n' \subset \Sigma_n$ be the maximal interval with $f^{\sigma(n)+(n-1)l+m_1}(\Sigma_n') \subset P_{m+1}' \cup P_m'$. Let $S_n(c_i) = \Sigma_n'$ and $s_i(n) = \sigma(n) + (n-1)l + m_1 + (m-1)n_1$. Note that $P_2' \subset f^{s_i(n)}(S_n(c_i))$.

Finally, (33) is easily seen to hold from the construction of $S_n(c_i)$. ∎

# 6  Lebesgue measure of limit sets

For a measurable set $\mathfrak{X}$, we denote by $|\mathfrak{X}|$ its Lebesgue measure. In this section we prove the following theorem.

**Theorem 6.1** *Let $f \in \mathfrak{D}$ and let $z \in I$. If $\omega(z)$ contains no intervals, then $|\omega(z)| = 0$.*

Such a theorem was proved for $C^2$ multimodal maps in [Lyu,1991], [Var,1996]. The proof is subdivided into several propositions, treating different kinds of limit sets. The $\omega$-limit set of a recurrent, nonperiodic, turning point is called a solenoid if it is obtained from infinitely many renormalizations. So, if $I_n$ is a sequence of decreasing intervals containing a turning point $c$ satisfying $f^{q(n)}(I_n) \subset I_n$ with $q_n \to \infty$ as $n \to \infty$ and $f^i(I_n) \cap f^j(I_n) = \emptyset$ for $0 \leq i < j < q(n)$, then $\omega(c) = \cap_{n \in \mathbb{N}} \cup_{0 \leq j < q_n} f^j(I_n)$ is a solenoid.

**Proposition 6.2** *Let $f \in \mathfrak{C}$. Let $z \in I$ satisfy $\mathfrak{T} \cap \omega(z) = \emptyset$, where $\mathfrak{T}$ denotes the set of turning points of $f$. Then $|\omega(z)| = 0$.*

PROOF. By extending $f$ to a larger interval and altering $f$ near the set $\mathfrak{T}$ of turning points, we may assume that $f(\partial I) \subset \partial I$ and $f(\mathfrak{T}) \subset \partial I$. This doesn't alter $\omega(z)$. Assume $\omega(z)$ has positive measure. Let $x \in \omega(z)$ be a point that is not eventually periodic. This excludes at most countably many points. We may thus take $x$ to be a density point of $\omega(z)$. Let $(y, y') \in I$ and $c \in \mathfrak{T}$ be such that $y \in \omega(x)$, $c \in (y, y')$ and $\omega(x) \cap (y, y') = \emptyset$. On a neighborhood $\mathcal{U}$ of $\{y\} \cup \{y'\}$ let a function $\tau$ be defined by $f(\tau(q)) = f(q)$ and $\tau(q) \neq q$. Note that $f(y) = y'$.

First assume that at most finitely many points in $\mathcal{O}(x)$ are contained in $(y, \tau(y))$. Replacing $x$ by an iterate, we may assume $\mathcal{O}(x) \cap (y, \tau(y)) = \emptyset$. Now let $f^{n(i)}(x)$ be the sequence of closest returns to $(y, \tau(y))$; $n(0)$ is the minimal integer with $f^{n(0)}(x) \in \mathcal{U}$ and $n(i+1)$ is the minimal integer with $f^{n(i+1)}(x) \in U_{f^{n(i)}(x)}$. Define

$$V_k = \{q \in \mathcal{U}, \quad f^i(q) \notin (q, \tau(q)), 0 < i < k, f^k(q) \in (q, \tau(q))\}, \tag{34}$$

see (9). Let $T_{n(i+1)-n(i)}$ be the component of $V_{n(i+1)-n(i)}$ containing $f^{n(i)}(x)$. Write $T_{n(i+1)-n(i)} = (a, b)$ with $b \in (y, a)$. Reasoning as in the proof of lemma 3.4 and noting that $f(\mathfrak{T}) \subset \partial I$, we get that either $f^{n(i+1)-n(i)}(a) = a$ and $f^{n(i+1)-n(i)}(b) = \tau(b)$ or $f^{n(i+1)-n(i)}(a) = \tau(a)$ and $f^{n(i+1)-n(i)}(b) = b$. Let $H_i = f^{n(i+1)-n(i)}(T_{n(i+1)-n(i)})$. Note that

$$(y, \tau(y)) \subset H_i. \tag{35}$$

Let $J_i$ be the maximal interval containing $x$ on which $f^{n(i)}(x)$ is monotone and $f^{n(i)}(J_i) \subset T_{n(i+1)-n(i)}$. We claim that

- $f^{n(i)}(J_i) = T_{n(i+1)-n(i)}$.

- $f^l(J_i) \cap f^k(J_i) = \emptyset$ for $0 \leq l < k < n(i)$.

The first item holds, since otherwise $f^j(p) = c$ for some $p \in J_i$, $j < n(i)$. Because $f^j(x) \notin U_a$, $f^j(J_i)$ contains either $T_{n(i+1)-n(i)}$ or $\tau(T_{n(i+1)-n(i)})$. So $f^{n(i)}(J_i)$ contains $f^{n(i)-j}(a)$ where $a \in \partial T_{n(i+1)-n(i)}$. It is however easily seen that $\mathcal{O}(a) \cap T_{n(i+1)-n(i)} = \emptyset$, a contradiction. To see that the intervals $f^j(J_i)$, $0 \leq j < n(i)$, are pairwise disjoint, suppose $f^l(J_i) \cap f^k(J_i) \neq \emptyset$ for some $0 \leq l < k < n(i)$. Then $f^{n(i)}(J_i) \cap f^{n(i)-k+l}(J_i) = T_{n(i+1)-n(i)} \cap f^{n(i)-k+l}(J_i) \neq \emptyset$. Since, by minimality of $n(i+1)$, $f^{n(i)-k+l}(J_i)$ can impossibly be contained in $T_{n(i+1)-n(i)}$, there exists $z \in J_i$ with $f^{n(i)-k+l}(z) = a \in \partial J_i$. This again contradicts $\mathcal{O}(a) \cap T_{n(i+1)-n(i)} = \emptyset$.

By lemma 2.3, the distortion of $f^{n(i+1)}$ on $J_i$ is bounded by some constant $D > 0$. From the forward invariance of $\omega(z)$ we get $f^{n(i+1)}(\omega(z) \cap J_i) \subset \omega(z) \cap H_i$. Hence, for some $q_i \in J_i$,

$$\begin{aligned}
\frac{|\omega(z) \cap H_i|}{|H_i|} &\geq \frac{|f^{n(i+1)}(\omega(z) \cap J_i)|}{|f^{n(i+1)}(J_i)|} \\
&= 1 - \frac{|f^{n(i+1)}(J_i \setminus (\omega(z) \cap J_i))|}{|f^{n(i+1)}(J_i)|} \\
&= 1 - \frac{\int_{J_i \setminus (\omega(z) \cap J_i)} |Df^{n(i+1)}(t)| dt}{|Df^{n(i+1)}(q_i)||J_i|} \\
&\geq 1 - D \frac{|J_i \setminus (J_i \cap \omega(z))|}{|J_i|},
\end{aligned}$$

which goes to 1 as $i \to \infty$ because $x$ is a density point of $\omega(z)$ and $|J_i| \to 0$, $i \to \infty$. So, because $(y, \tau(y)) \subset V_i$, $|\omega(z) \cap (y, \tau(y))| = |(y, \tau(y))|$. Since $\omega(z)$ is a closed set, we get $(y, \tau(y)) \subset \omega(z)$, contradiction. The case where infinitely many iterates $f^i(x)$ are contained in $(y, \tau(y))$ is treated similarly, compare the proof of lemma 3.3. ∎

Before going on studying the measure of $\omega$-limit sets of recurrent turing points, we introduce some notions and notations. Near a turning point $c$, one can define a function $\tau$ by demanding $f(c) = c$ and $f(\tau(y)) = f(y)$ with $\tau(y) \neq y$ is $y \neq c$. Let $\mathcal{U}_\tau$ be the union of the neighborhoods of turning points on which such a function is defined. For $x \in \mathcal{U}_\tau$, write $U_x = (x, \tau(x))$. A point $x \in \mathcal{U}_\tau$ is called *nice* if $\mathcal{O}(x) \cap U_x = \emptyset$. An interval of the form $U_x$ is called *symmetric*.

The following lemma gives the device by which one can prove that the $\omega$-limit set of a recurrent turning point, if it contains no intervals, has zero Lebesgue measure.

**Lemma 6.3** *Let $f \in \mathfrak{C}$ and let $c$ be a turning point of $f$. If there exist $\lambda > 1$ and symmetric neighborhoods $P_n \subset Q_n$ of $c$ with $|Q_n| \to 0$ as $n \to 0$, $\partial Q_n$ consisting of nice points, $\omega(c) \cap Q_n \subset P_n$, and $|Q_n|/|P_n| \geq \lambda$, then $|\omega(c)| = 0$.*

PROOF. Let $\mathfrak{C} = \{x, \ c \in \omega(x)\}$. By theorem 6.1, $\omega(c) \cap (I \backslash \mathfrak{C})$ has zero measure. Take $x \in \mathfrak{C} \cap \omega(c)$ and let $h(n)$ be the minimal integer with $f^{h(n)}(x) \in Q_n$. Let $T_n$ be the maximal interval containing $x$ so that $f^{h(n)}(T_n) \subset Q_n$. We claim that

- $f^i(T_n) \cap f^j(T_n) = \emptyset$ for $0 \leq i < j < h(n)$.

Indeed, if $f^i(T_n) \cap f^j(T_n) \neq \emptyset$ for some $0 \leq i < j < h(n)$, then $f^{h(n)}(T_n) \cap f^{h(n)-j+i}(T_n) \neq \emptyset$. Because $f^{h(n)-j+i}(T_n) \subset Q_n$ is not possible by minimality of $h(n)$, there exists $z \in T_n$ with $f^{h(n)-j+i}(z) \in \partial Q_n$. But then $f^{h(n)}(z)$ cannot be in $Q_n$ since $\partial Q_n$ consists of nice points, contradiction.

Let $S_n \subset T_n$ satisfy $f^{h(n)}(S_n) \subset P_n$. Note that

$$\omega(c) \cap T_n \ \subset \ S_n. \tag{36}$$

By lemma 2.3, $f^{h(n)}$ has bounded distortion on $T_n$ with a bound independent of $n$. Therefore, using $\partial f^{h(n)}(T_n) \subset \partial Q_n$ and $|Q_n|/|P_n| \geq \lambda$, there exists $\tilde{\lambda} > 1$ so that

$$|T_n|/|S_n| \ \geq \ \tilde{\lambda} \tag{37}$$

for all $n$. The contraction principle lemma 2.1 yields $|T_n| \to 0$ as $n \to \infty$. Therefore, (36) and (37) imply that $x$ is not a density point of $\omega(c)$. So almost no point in $\omega(c)$ is a density point of $\omega(c)$. By the Lebesgue density theorem, $\omega(c)$ has zero measure. ∎

**Proposition 6.4** *Let $f \in \mathfrak{C}$. Let $c$ be a turning point of $f$.*
*If $\omega(c)$ is not minimal and does not contain an interval, then $|\omega(c)| = 0$.*

PROOF. Because $\omega(c)$ is not minimal, there exists $y \in \omega(c)$ with $c \notin \omega(y)$. We may assume that $\omega(y)$ is a minimal set.

First suppose $y$ is periodic. Write $k$ for the minimal period of $y$. If $Df^k(y) < 0$, let $l = 2k$, otherwise let $l = k$. Let $I_1 = (y, y_1)$ be an interval on which $f^l$ is monotone, with $y_1$ an eventually periodic point and $\mathcal{O}(c)$ accumulates in $I_1$ to $y$. Let $I_n$ be the maximal subinterval of $I_1$ with $f^{(n-1)l}(I_n) = I_1$. Write $J_n = I_n \backslash I_{n-1}$. Note that $J_n$ is a fundamental domain; $f^l(\text{int } J_n) \cap \text{int } J_n = \emptyset$ and $f^l(\partial I_n) \cap \partial I_n \neq \emptyset$.

Let $s(n)$ be the minimal integer with $f^{s(n)}(f(c)) \subset J_n$. Since $\mathcal{O}(c)$ accumulates in $I_1$ on $y$, $s(n)$ is well defined for all $n$. Let $S_n$ be the maximal interval containing $f(c)$ so that $f^{s(n)}(S_n) \subset J_n$. We claim that for $n$ sufficiently large,

- $f^i(S_n) \cap f^j(S_n) = \emptyset$ for $0 \leq i < j \leq s(n)$.

Indeed, if $f^i(S_n) \cap f^j(S_n) = \emptyset$ for some $0 \leq i < j \leq s(n)$, then $f^{s(n)}(S_n) \cap f^{s(n)-j+i}(S_n) \neq \emptyset$. Now $f^{s(n)-j+i}(S_n) \subset J_n$ is impossible by minimality of $s(n)$. So some $p \in S_n$ is mapped by $f^{s(n)-j+i}$ to $z \in \partial J_n$. This contradicts $\mathcal{O}(z) \cap J_n = \emptyset$ for $n$ high which follows from the fact that $y_1$ is eventually periodic.

Applying lemma 2.3 to $f^{s(n)}|_{S_n}$ and $f^{ln}|_{J_n}$, one sees that $f^{s(n)+ln}$ has uniformly bounded distortion on $S_n$.

By taking $y_1$ with $\omega(y_1) \cap \omega(c) = \emptyset$, we have $\partial J_1 \cap \omega(c) = \emptyset$. Let $T \subset J_1$ be so that $\omega(c) \cap J_1 \subset T$, let $T_n \subset S_n$ be so that $f^{s(n)+(n-1)l}(T_n) = T$. Because $f^{s(n)+(n-1)l}$ has bounded distortion on $S_n$ and $\partial f^{s(n)+(n-1)l} S_n \subset \partial J_1$, there exists $\nu > 0$ with

$$|L_n|/|T_n|, |R_n|/|T_n| \geq \nu, \tag{38}$$

where $L_n, R_n$ are the components of $S_n \setminus T_n$. Let $P_n, Q_n$ be the maximal intervals containing $c$ with $f(P_n) \subset T_n$, $f(Q_n) \subset S_n$ respectively. By (38), $|Q_n|/|P_n| \geq \lambda$ for some $\lambda > 1$. An application of lemma 6.3 yields $|\omega(c)| = 0$.

The case where $y$ is an infinite minimal set proceeds similarly. We may assume $y$ is a nice point. Let $f^{n(i)}(f(y))$ be the sequence of closest returns to $c$; $n(0)$ is the minimal integer with $f^{n(0)}(f(y))$ in the connected component of $\mathcal{U}_\tau$ that contains $c$ and $n(i+1)$ is the minimal integer with $f^{n(i+1)}(f(y)) \in U_{f^{n(i)}(y)}$. Let $V_k$ be defined as in (34) and let $T_{n(i+1)-n(i)}$ be the component of $V_{n(i+1)-n(i)}$ containing $f^{n(i)}(f(y))$. Write $T_{n(i+1)-n(i)} = (a,b)$ with $b$ closer to $c$. Let $J_n$ be the maximal interval containing $f(y)$ with $f^{n(i)}(f(y)) \subset T_{n(i+1)-n(i)}$. From the fact that a boundary point of $T_{n(i+1)-n(i)}$ or its image under $f$ is periodic, see the proof of lemma 3.4, it is easily seen that

$$\mathcal{O}(\partial T_{n(i+1)-n(i)}) \cap T_{n(i+1)-n(i)} = \emptyset. \tag{39}$$

We claim that

- $f^k(J_n) \cap f^l(J_n) = \emptyset$ for $0 \leq k < l < n(i)$.

Otherwise $f^{n(i)}(J_n) \cap f^{n(i)-l+k}(J_n) = T_{n(i+1)-n(i)} \cap f^{n(i)-l+k}(J_n) \neq \emptyset$. Since, by minimality of $n(i+1) - n(i)$, $f^{n(i)-l+k}(J_n) \subset T_{n(i+1)-n(i)}$ is not possible, there is $z \in J_n$ with $f^{n(i)-l+k}(z) \in \partial T_{n(i+1)-n(i)}$. Now (39) yields a contradiction. The rest proceeds as above. ∎

**Proposition 6.5** *Let $f \in \mathfrak{C}$. Let $c$ be a turning point of $f$.*

*If $c$ is a solenoid, then $|\omega(c)| = 0$.*

PROOF. This proposition is proved in [AlsLopSno,1995]. For completeness, we give a short alternative argument here. The map $f$ restricted to $\omega(c)$ is injective. Therefore, for each $y \in \mathcal{U}_\tau$, $|\omega(c) \cap U_y|$ is smaller than the maximum of $|(y,c)|$ and $|(c,\tau(y))|$. So

$$\limsup_{y \to c} |\omega(c) \cap U_y|/|U_y| \leq \lambda \tag{40}$$

for some $\lambda < 1$. Take $x \in \omega(c)$. If $|\omega(c)| > 0$ we may assume $x$ is a density point of $\omega(c)$. Let $c_1, \ldots, c_s$ with $c_s = c$ denote the turning points in $\omega(c)$. Since $f$ is infinitely renormalizable, there is a decreasing sequence of neighborhoods $P_n(c_i)$ of $c_i$ with $f^{q_n}(P_n(c_i)) \subset P_n(c_i)$ and $f^k(P_n(c_i)) \cap f^l(P_n(c_i)) = \emptyset$, $0 \leq k < l < q(n)$, for some $q(n) \to \infty$ as $n \to \infty$. Let $k(n) > 0$ be the minimal number such that $f^{k(n)}(x) \in P_n(c_i)$ for some $c_i$. Let $I_n$ be the maximal interval so that $f^{k(n)}$ maps $I_n$ homeomorphically onto $P_n(c_i)$. By lemma 2.3, there exists $D > 0$ so that the distortion of $f^{k_n}|_{I_n}$ is bounded by $D$. Compute

$$\frac{|P_n(c_i)\setminus(P_n(c_i) \cap \omega(c))|}{|P_n(c_i)|} \leq \frac{|f^{k_n}(I_n\setminus(I_n \cap \omega(c)))|}{|f^{k_n}(I_n)|} \leq D\frac{|I_n\setminus(I_n \cap \omega(c))|}{|I_n|}.$$

Since $x$ is a density point of $\omega(c)$ and $|I_n| \to 0$ as $n \to \infty$, $|I_n\setminus(I_n \cap \omega(c))|/|I_n| \to 0$ as $n \to \infty$. Thus $|P_n(c_i) \cap \omega(c)|/|P_n(c_i)| \to 1$ as $n \to \infty$. This contradicts (40). ∎

The following proposition is proved in [Var,1996] for $C^2$ multimodal maps. His proof also works for $f \in \mathfrak{D}$. We present a proof for the smaller class $\mathfrak{E}$ of multimodal maps $f \in \mathfrak{D}$ for which $|Df|$ is continuous and refer to [Var,1996] for the proof for $f \in \mathfrak{D}$.

**Proposition 6.6** *Let $f \in \mathfrak{D}$ and let $c$ be a turning point of $f$. If $\omega(c)$ is a minimal Cantor set which is not a solenoid, then $|\omega(c)| = 0$.*

A main ingredient of the proof of proposition 6.6 is lemma 6.7 below, for which we first introduce some notation.

For $x \in \mathcal{U}_\tau$, let $D_x = \{y \in I, \exists k > 0 \text{ with } f^k(y) \in U_x\}$. The minimal number $k$ with $f^k(y) \in U_x$ for $y \in D_x$ is called the *transfer time* of $y$. Let $R_x : D_x \to U_x$ be the Poincaré map; $R_x(y) = f^k(y)$ where $k$ is the transfer time of $y$. If $x$ is a nice point then the transfer time is constant on each connected component of $D_x$: a boundary point of a maximal subinterval of $D_x$ on which the transfer time is constant, is mapped to $\partial U_x$ for some iterate and thus is also in $\partial D_x$.

Let $z$ be some nice point and write $S_z$ for the connected component of $D_z$ that contains $f(c)$. Let $x = \psi(z)$ be defined by $U_x = f^{-1}(S_z)$.

**Lemma 6.7** *Let $f \in \mathfrak{E}$. Suppose $f$ is not renormalizable and $c$ is recurrent. There are symmetric neighborhoods $U_{u(n)} \subset U_{v(n)}$ of $c$, where $v(n)$ a nice point and $u(n) = \psi(v(n))$ as defined above, with the following properties.*

1. $|U_{v(n)}| \to 0$ as $n \to \infty$.

2. $f^i(u(n)) \notin U_{v(n)}$ for $i > 0$.

3. There exists $D > 0$ so that the distortion of $R_{u(n)}$ is bounded on each connected component of $D_{u(n)}$.

4. For some $\lambda < 1$, $|U_{u(n)}|/|U_{v(n)}| \leq \lambda$.

PROOF. Let $v(n)$ be a nice point and let $u(n) = \psi(v(n))$. For the second item, first note that $f^k(u(n)) = v(n)$, where $R_{u(n)}|_{U_n} = f^k$, so that $f^i(u(n)) \notin U_{v(n)}$ if $i \geq k$. If $f^i(u(n)) \in U_{v(n)}$ for some $0 < i < k$, then since $f^i(u(n)) \notin U_{v(n)}$, there is $a \in U_{u(n)}$ with $f^i(a) \in \{v(n), \tau(v(n))\}$. But then $f^k(a) \in U_{v(n)}$ is impossible since $v(n)$ is nice. For the third item, showing bounded distortion of $R_{u(n)}$ on connected components of $D_{u(n)}$, let $J$ be a connected component of $D_{u(n)}$ and let $k$ denote the transfer time on $J$. Then $f^i(J) \cap f^j(J) = \emptyset$ for $0 \leq i < j < k$. To see this, note that $f^i(J) \cap J \neq \emptyset$ for some $0 < i < k$ would imply the existence of $y \in J$ with $f^i(y) = a \in \partial J$. Then $f^k(y) = f^{k-i}(a) \notin U_{u(n)}$, a contradiction. Lemma 2.3 shows that, for some $D > 0$, the distortion of $R_x$ on each connected component of $D_{u(n)}$ is bounded by $D$.

For the last item, we must further restrict the choice of $v(n)$. Let $c_{q(n)} = f^{q(n)}(c)$ be the sequence of nearest returns to $c$. Denote by $\mathfrak{N}$ the set of nice points; observe that $c$ is an accumulation point of nice points. There thus are infinitely many integers $n$ with $(U_{c_{q(n-1)}} \setminus U_{c_{q(n)}}) \cap \mathfrak{N} \neq \emptyset$. For those $n$, let

$$z(n) = \sup\{y < c, \ y \in (U_{c_{q(n-1)}} \setminus U_{c_{q(n)}}) \cap \mathfrak{N}\},$$
$$x(n) = \psi(z(n)).$$

By definition of $z(n)$, we have $x(n) \in U_{c_{q(n)}}$ and $R_{x(n)}(c) = c_{q(n)} \in U_{z(n)} \setminus U_{x(n)}$.

$R_{x(n)}$ *has a fixed point in* $U_{x(n)}$. Let $p_{q(n)} \in U_{x(n)}$ be the periodic point of minimal period $q(n)$ and and $Df^{q(n)}(p_{q(n)}) > 0$. Since $f$ is not renormalizable, there is $r(n) \in (p_{q(n)}, c)$ with $f^{q(n)}(r(n)) = \tau(p_{q(n)})$. Since the distortion of $R_{z(n)}$ on $U_{x(n)}$ is uniformly bounded and $R_{z(n)}(c) \in U_{z(n)} \setminus U_{x(n)}$, there exists $\nu > 1$ so that if $|U_{z(n)}|/|U_{x(n)}| \leq \nu$, then $|U_{p_{q(n)}}|/|U_{r(n)}| \geq \nu$. Therefore we can always choose neighborhoods $U_{u(n)} \subset U_{v(n)}$ as in the statement of the lemma.

$R_{x(n)}$ *has no fixed point in* $U_{x(n)}$. We claim the existence of $\nu > 1$ so that if $|U_{z(n)}|/|R_{x(n)}(U_{x(n)})| \leq \nu$, then $c$ is in the basin of attraction of a periodic attractor. If $|U_{z(n)}|/|R_{x(n)}(U_{x(n)})| \leq \nu$, then $|DR_{x(n)}|_{U_{x(n)}}| \leq D(\nu - 1)$ for some constant $D$, since $R_{x(n)}$ has uniformly bounded distortion on $U_{x(n)}$. Applying lemma 2.5 one easily obtains the claim.

So we may assume $|U_{z(n)}|/|R_{x(n)}(U_{x(n)})| \leq \nu$. Let $\bar{\nu} < \nu$. If $|U_{z(n)}|/|U_{x(n)}| \geq \bar{\nu}$ let $v(n) = z(n)$ and $u(n) = x(n)$. Otherwise let $v(n)$ be so that $U_{v(n)}$ is the maximal interval in $U_{x(n)}$ with $f^{q(n)}(U_{x(n)} \setminus U_{v(n)}) \subset U_{z(n)} \setminus U_{x(n)}$. Let $r(n) = \Psi(q(n))$. Because $R_{x(n)}$ has uniformly bounded distortion on $U_{x(n)}$, we have $|U_{q(n)}|/|U_{r(n)}| \geq \nu'$ for some $\nu' > 1$. ∎

PROOF OF PROPOSITION 6.6 FOR $f \in \mathfrak{E}$. Take symmetric neighborhoods $U_{u(n)} \subset U_{v(n)}$ as in lemma 6.7. Since $\omega(c)$ is minimal, $\omega(c) \subset D_{u(n)}$. Compactness of $\omega(c)$ implies that $\omega(c)$ is contained in a finite number of connected components of $D_{u(n)}$. Let $I_n$ be the connected component with maximal transfer time $k_n$. Denote by $J_n \supset I_n$ the maximal interval with $f^{k_n}(J_n) \subset U_{v(n)}$. We claim that

$$\omega(c) \cap (J_n \setminus I_n) = \emptyset. \tag{41}$$

If this were not so, there would be an $y \in J_n \setminus I_n$ with $f^i(y) \in U_{u(n)}$, $i < k_n$ (by maximality of $k_n$). Then $f^i(I_n) \cap U_{u(n)} = \emptyset$ and $f^i(y) \in U_{u(n)}$ implies that $f^i(z) \in \{u(n), \tau(u(n))\}$ for some $z \in J_n \setminus I_n$. Then $f^{k_n}(z) = f^{k_n - i}(u(n))$. By lemma 6.7, $f^j(u(n)) \notin U_{v(n)}$ for all $j > 0$. So $f^{k_n}(z)$ can not be in $U_{u(n)}$. This contradiction shows (41).

Let $s_n$ be the minimal integer with $f^{s_n}(f(c)) \in I_n$. As shown before there is an interval $S_n$ containing $f(c)$ so that $f^{s_n}(S_n) \subset J_n$. We claim that

- $f^k(S_n) \cap f^l(S_n) = \emptyset$ for $0 \leq l < k < s_n$.

If $f^k(S_n) \cap f^l(S_n) \neq \emptyset$ for some $0 \leq l < k < s_n$, then $f^{s_n}(S_n) \cap f^{s_n-l+k}(S_n) \neq \emptyset$. Write $h = s_n - l + k$ and observe that $f^h(S_n)$ can not be contained in $J_n = f^{s_n}(S_n)$ by minimality of $s_n$. So there exists $z \in J_n$ with $f^h(z) = a \in \partial J_n$. Then $f^{h+k_n}(z) \in \{v(n), \tau(v(n))\}$. So $f^{s_n+k_n}(z) = f^{s_n+k_n-h}(a) = f^{s_n-j}(v(n))$ can not be in $U_{v(n)}$ since $v(n)$ is nice, a contradiction.

By lemma 2.3, $f^{s_n}$ has bounded distortion on $S_n$ where the bound does not depend on $n$. By lemma 6.7, $|U_{u(n)}|/|U_{v(n)}| \leq \lambda$ for some $\lambda < 1$ not depending on $n$ and the map $f^{k_n}$ has bounded distortion on $I_n$. Let $\tilde{J}_n \subset J_n$ be the maximal interval with $\sum_{i=0}^{k_n-1} |f^i(\tilde{L}_n)|$ and $\sum_{i=0}^{k_n-1} |f^i(\tilde{R}_n)|$ bounded by $\sum_{i=0}^{k_n-1} |f^i(I_n)|$, where $\tilde{L}_n, \tilde{R}_n$ are the components of $\tilde{J}_n \setminus I_n$. By lemma 2.4, $f^{k_n}$ has uniformly bounded distortion on $\tilde{J}_n$. For some $\nu > 0$, both $f^{k_n}(\tilde{L}_n)$ and $f^{k_n}(\tilde{R}_n)$ have size at least $\nu |U_{u(n)}|$. Let $P'_n \subset Q'_n$ be intervals containing $f(c)$ with $f^{k_n+s_n}(P'_n) = U_{u(n)}$, $f^{k_n+s_n}(Q'_n) = U_{v(n)}$. Replacing $Q'_n$ by the smaller interval $\tilde{Q}_n$ satisfying $f^{s_n}(\tilde{Q}_n) = \tilde{J}_n$ we can ensure that $f^{k_n+s_n}$ has bounded distortion on $\tilde{Q}_n$ and both components of $\tilde{Q}_n \setminus P'_n$ have size at least $\nu |P'_n|$ for some $\nu > 0$. So also

$$|L'_n|/|P'_n|, |R'_n|/|P'_n| \geq \nu,$$

where $L'_n, R'_n$ are the components of $Q'_n \setminus P'_n$. Applying lemma 6.3 with $P_n = f^{-1}(P'_n)$ and $Q_n = f^{-1}(Q'_n)$ (well defined since $f(c) \in P'_n$) yields $|\omega(c)| = 0$. ∎